\documentclass[12pt, final]{amsproc}
\usepackage{amssymb}
\usepackage[mathscr]{eucal}
\usepackage{amscd}
\usepackage{amsmath}
\usepackage{cite}
\usepackage{setspace}
\usepackage{tikz}
\usetikzlibrary{shapes,arrows,fit,calc,positioning}
\tikzset{box/.style={draw, rectangle, thick, text centered, minimum height=3em}}
\tikzset{line/.style={draw, thick, -latex'}}
\def\Ker{\operatorname{Ker}}

\theoremstyle{plain}
\newtheorem{theorem}{Theorem}

\newtheorem*{psp}{The Principle of Small Perturbations}

\newtheorem{example}{Example}
\theoremstyle{remark}

\title{A remark on the distortion of twisted sums}
\address{Departamento de Matem\'aticas, Universidad de Extremadura, Avenida de Elvas s/n, 06011 Badajoz, Spain.}
\email{jesus@unex.es}
\author[J. Su\'arez]{Jes\'us Su\'arez}
\thanks{The author was supported by no institution.}
\begin{document}
\begin{abstract} Odell and Schlumprecht solved the distortion problem by proving that the classic sequence spaces $\ell_p$ for $1<p<\infty$ admit an inevitable biorthogonal system. In particular, these spaces are arbitrarily distortable. Later, Maurey extended the result to asymptotic $\ell_p$-spaces while Tomczak-Jaegermann did likewise for the Schatten classes. We observe that the Kalton-Peck spaces $Z_p$ for $1<p<\infty$ admit an inevitable biorthogonal system. Therefore, we may add this classic family to the list of arbitrarily distortable spaces. 
\end{abstract}
\maketitle
\section{Introduction}% This means roughly that every Lipschitz function on the unit sphere of these two spaces must stabilize, this is, the oscillation is as small as we wish if we restrict ourselves to the unit sphere of a suitable subspace. Taking as a Lipschitz function any equivalent norm, we have the formal definition.
As it is well known, it was James who introduced us to the topic of the distortion. He proved that the spaces $\ell_1$ and $c_0$ are not distortable \cite{Ja}. Recall that given $\lambda >1$, a norm $\|\cdot\|$ on an infinite-dimensional Banach space $X$ is said to be $\lambda$-distortable if there is an equivalent norm $|||\cdot|||$ on $X$ such that for every infinite-dimensional subspace $Y$ of $X$ $$\sup\{|||x|||/|||y|||: x,y\in Y, \|x\|=\|y\|=1\}\geq \lambda.$$
A given norm is said to be distortable if it is $\lambda$-distortable for some $\lambda >1$ and arbitrarily distortable if it is $\lambda$-distortable for all $\lambda>1$. 
We invite the reader to the survey of Odell and Schlumprecht \cite{OS2} for a comprehensive study on distortion and its development. We just focus on some points to introduce our results, starting with Schlumprecht's construction in \cite{Sch} of the first example of a Banach space which is arbitrarily distortable and that we denote by $\mathcal S$. This follows nowadays from the fact that $\mathcal S$ admits an inevitable biorthogonal system. Here, we are replacing the term asymptotic by inevitable as suggested by Maurey in \cite{Mau}. Recall that a Banach space $X$ is said to admit an inevitable biorthogonal system with constant $0<\delta<1/2$ if there are two sequences $\{A_j\}_{j=1}^\infty$ and $\{A_j^*\}_{j=1}^\infty$ of subsets of the unit sphere $S(X)$ of $X$ and of the unit ball $B(X^*)$ of $X^*$, respectively, such that
\begin{enumerate}
\item $d(H,A_j)=0$ for every infinite dimensional subspace $H$ of $X$ and every $j\in \mathbb N$.
\item Given $j\in \mathbb N$ and $x\in A_j$, there is $x^*\in A_j^*$ so that $\langle x^*,x\rangle \geq 1-\delta$.
\item If $j\neq k$, then $\langle x^*, x\rangle \leq \delta$ whenever $x\in A_j$ and $x^*\in A_k^*$.
\end{enumerate} 
%The term inevitable is justified by condition (1) above since every $A_j$ must essentially  meet the unit sphere of every infinite-dimensional subspace of $X$. The other two conditions reflect obviously the meaning of biorthogonal. 
The condition (1) above is usually termed as the set $A_j$ is inevitable for $j\in \mathbb N$. The way to construct such a system in $\mathcal S$ deals with $\ell_1$-averages and rapidly increasing sequences sequences, RIS sequences in short, which are tools introduced by Gowers and Maurey in the solution to the unconditional basic sequence problem in \cite{GM}. Once we have such a system for a fixed $\delta$, the equivalent norm
$$\varepsilon\|x\|+\sup\{|x^*(x)|:x^*\in A_1^*\},$$
is a $(1+\varepsilon-\delta)/(\varepsilon+\delta)$-distortion, see \cite[Proposition 13.28.]{BenLin} for the details.

A natural question remained opened for a long time: Whether or not the classical sequence spaces $\ell_p$ for $1<p<\infty$ are also arbitrarily distortable. It was a surprise when Odell and Schlumprecht in \cite{OS} were able to transfer these inevitable sets from $\mathcal S$ to $\ell_p$ with $1<p<\infty$ solving the ``the distortion problem". Later, Maurey extended the result to asymptotic $\ell_p$-spaces with an unconditional basis \cite{Mau2}. These are spaces such that any $n$-blocks of the basis starting  far enough are uniformly equivalent to the unit vector basis of $\ell_p$. Modelled over Maureys' proof, it was shown by Tomczak-Jaegermann in \cite{TJ} that the Schatten classes $C_p$ also admit such system of inevitable sets for $1<p<\infty$.  

A third class of spaces which fits with the previous idea is that of the twisted sums of $\ell_p$. These are Banach spaces $Z$ so that they contain a copy of $\ell_p$ for which the natural quotient map $Z/\ell_p$ is again isomorphic to $\ell_p$. The main examples of this class are the so-called Kalton-Peck spaces $Z_p$ for $1<p<\infty$ which were introduced in the seminal paper \cite{KaPe}. The main purpose of this note is to extend the inevitable biorthogonal system from $\ell_p$ to a large class of twisted sums of $\ell_p$ that includes the Kalton-Peck spaces. So we add new examples to the works of Maurey \cite{Mau2} and Tomczak-Jaegermann \cite{TJ} on families of arbitrarily distortable Banach spaces.

% Our contribution in this short note is to complete the works of Maurey \cite{Mau2} and Tomczak-Jaegermann \cite{TJ} proving the equivalent results for the Kalton-Peck spaces.
\section{Main}
We shall use the terminology commonly used in Banach space theory as it appears in the book of Albiac and Kalton \cite{AK} where all unexplained notation may be found.

Let us briefly recall the formal construction of a twisted sum of $\ell_p$. First, an homogeneous map $\Omega:\ell_p\longrightarrow \omega$ is called a centralizer if there is a constant $C>0$ so that: 
\begin{equation*}\label{centralizer}
\| \Omega(ax)-a\Omega(x)\|_{\ell_p}\leq C\|a\|_{\infty}\|x\|_{\ell_p},\;\;\:\;a\in \ell_{\infty},x\in \ell_p,
\end{equation*}
where $\omega$ denotes the vector space of all complex scalar sequences. By homogeneous we just mean that $\Omega$ must satisfy $\Omega(\lambda x)=\lambda\Omega(x)$ for every $\lambda \in \mathbb C$ and $x\in \ell_p$. Once we have a centralizer $\Omega$, we pick the pairs $(x,y)\in \omega \times \omega$ for which the quasi-norm $$\|x-\Omega(y)\|_{\ell_p}+\|y\|_{\ell_p}$$
is finite. The quasi-norm above is equivalent to a norm for $1<p<\infty$ by a result of Kalton \cite[Theorem 2.6]{ka3}. This set of couples, say $Z$, is a twisted sum of $\ell_p$. Indeed, it is easy to see that $j(x)=(x,0)$ is an embedding of $\ell_p$ into $Z$ for which the quotient map is identified as $q(x,y)=y$ from $Z$ onto $\ell_p$. The most important centralizer is given by the so-called Kalton-Peck map 
\begin{equation}\label{kpmap}
\Omega_p(x)=\sum_{j=1}^\infty x_j\cdot \log \frac{|x_j|}{\|x\|_p}
\end{equation}
with the convention $0\cdot \log 0=0$. The Kalton-Peck space $Z_p$ is thus defined as the pairs $(x,y)\in \omega \times \omega$ for which the quasi-norm $$\|x-\Omega_p(y)\|_{\ell_p}+\|y\|_{\ell_p}$$
is finite. The details are to be found in the survey of Kalton and Montgomery-Smith \cite[Pages 1155--61]{KM}. It will be useful for us to recall the duality between a twisted sum of $\ell_p$, say $Z$, and its dual $Z^*$. The natural pairing for $(x,y)\in Z$ and $(a,b)\in Z^*$ is given by the rule $$\langle(x,y),(a,b)\rangle=\langle x,b\rangle+ \langle y, a\rangle.$$
The proof of this fact for the Kalton-Peck map is \cite[Theorem 5.1]{KaPe}. For an arbitrary centralizer and $p=2$, a proof can be found in \cite[Lemma 1]{JS2} but the proof easily extends to $1<p<\infty$.
The class of twisted sums we work with is isolated by the following rule.
\begin{psp} A twisted sum $Z$ of $\ell_p$ for $1<p<\infty$ satisfies P.S.P. if for every sequence $(w_j)_{j=1}^{\infty}$ in $Z$ which is equivalent to the unit vector basis of $\ell_p$ there is a subsequence $(w_{k_j})_{j=1}^{\infty}$ of $(w_j)_{j=1}^{\infty}$ and blocks $(y_j)_{j=1}^{\infty}$ in $\ell_p$ such that 
\begin{enumerate}
\item $\|(y_j,0)\|=1$ for every $j\in \mathbb N$.
\item $\|w_{k_j}-(y_j,0)  \|\leq 2^{-j}$ for every $j\in \mathbb N$.
\end{enumerate}
\end{psp}
The Principle of Small Perturbations, the P.S.P. in short, above is satisfied by the Kalton-Peck spaces as we will show in Example 1 in the next section. In the meantime, let us justify the importance of the P.S.P. for the distortion.
\begin{theorem}\label{main}
Let $Z$ be a twisted sum of $\ell_p$ with $1<p<\infty$. If $Z$ satisfies the Principle of Small Perturbations, then $Z$ admits an inevitable biorthogonal system for every $0<\delta<1/2$. In particular, $Z$ is arbitrarily distortable.
\end{theorem}
\begin{proof}
Let $\{A_j,A_j^*\}_{j=1}^{\infty}$ be an inevitable biorthogonal system of $\ell_p$ with constant $\delta>0$. We define sets in $Z$ for every $j\in \mathbb N$ as
$$\widehat{A}_j:=\{(x,0):x\in A_j\}.$$
We turn to $Z^*$ so let $\Omega$ be a centralizer for $Z^*$ and define the following sets in $Z^*$ for every $j\in \mathbb N$
$$\widehat{A}_j^*:=\{(\Omega(x),x):x\in A_j^*\}.$$
Observe that the quasinorm of an element $(\Omega(x),x)$ with $x\in A_j^*$ is given by $$\|\Omega(x)-\Omega(x)\|+\|x\|=\|x\|<\infty,$$
so the sets $\widehat{A}_j^*$ are well defined.
The proof consists in showing that the sets $\{\widehat{A}_j\}_{j=1}^\infty$ and $\{\widehat{A}_j^*\}_{j=1}^\infty$ defined above form an inevitable biorthogonal system for $Z$. We only need to check that the set $\widehat{A}_j$ is inevitable for each $j\in\mathbb N$. Indeed, the reimaining two claims of the definition are trivial using the natural pairing duality $$\langle(x,y),(a,b)\rangle=\langle x,b\rangle+ \langle y, a\rangle$$ and the fact that $\{A_j,A_j^*\}_{j=1}^{\infty}$ is biorthogonal. Thus, pick $j_0\in \mathbb N$ and let $W$ be any infinite dimensional subspace of $Z$. Our aim is to show that $d(W,\widehat{A}_{j_0})=0$. Thus, let us fix $\varepsilon>0$.
Let us recall first that every infinite-dimensional subspace $W$ of a twisted sum of $\ell_p$ contains an isomorphic copy of $\ell_p$: First, every infinite dimensional subspace of $\ell_p$ contains an isomorphic copy of $\ell_p$, see \cite[Proposition 2.2.1]{AK}. Now, if $W\cap \Ker q$ is infinite-dimensional, where $q(a,b)=b$ denotes the natural quotient map onto $\ell_p$, this is clear by the first argument. Otherwise, $q$ is an isomorphism on a finite codimensional subspace of $W$ and we are done again. This argument is just taken for the case $p=2$ from the proof of \cite[Theorem 16.16., Second paragraph]{BenLin}. By the P.S.P., we may find a subsequence $(w_{k_n})_{n=1}^{\infty}$ of $(w_n)_{n=1}^{\infty}$ and normalized blocks $(u_n)_{n=1}^{\infty}$ in $\ell_p$ so that $$\|w_{k_n}-(u_n,0)\|\leq \frac{\varepsilon}{2^n},\;\;\;n\in \mathbb N.$$
Define $U$ as the closed linear span of $(u_j)_{j=1}^{\infty}$ in $\ell_p$, usually termed as $U=[u_n]_{n=1}^{\infty}$ . Since $A_{j_0}$ is inevitable in $\ell_p$, given our $\varepsilon>0$ and $U$, there must be a vector $u\in U$ with $\|u\|_{\ell_p}\leq 1$ and $a\in A_{j_0}$ so that $$\|u-a\|_{\ell_p}\leq \varepsilon.$$
We now consider $Z$, so write $u=\sum_{n=1}^{\infty} \lambda_n u_n$ and thus let $w=\sum_{n=1}^{\infty} \lambda_n w_{k_n}\in W$. Observe that $w$ is a well defined element of $W$. Indeed, $\|u\|_{\ell_p}\leq 1$ and the normalized blocks $(u_n)_{n=1}^{\infty}$ are equivalent to the canonical $\ell_p$ basis. This means that $$\|u \|_{\ell_p}^p=\sum_{n=1}^{\infty} |\lambda_n|^p\|u_n\|_{\ell_p}^p=\sum_{n=1}^{\infty} |\lambda_n|^p\leq 1.$$ Since $(w_{k_n})_{n=1}^{\infty}$ is again equivalent to the canonical $\ell_p$-basis, a similar computation as before proves our claim. Then, we have
\begin{eqnarray*}
\|w-(a,0)\|_Z&\leq & \|w-(u,0)\|_Z+\|(u,0)-(a,0)\|_Z\\
&=&  \|w-(u,0)\|_Z+\|u-a\|_{\ell_p}\\
&\leq&  \|w-(u,0)\|_Z+\varepsilon\\
&\leq&  \sum_{n=1}^{\infty}|\lambda_n| \|w_{k_n}-(u_n,0)\|+\varepsilon\\
&\leq&  \sum_{n=1}^{\infty} \frac{\varepsilon}{2^n}+\varepsilon\\
&=& 2\varepsilon.
\end{eqnarray*}
Since obviously $(a,0)\in \widehat{A}_{j_0}$, the inequality above gives that $d(W,\widehat{A}_{j_0})\leq 2\varepsilon$. But $\varepsilon>0$ is arbitrary,  so $d(W,\widehat{A}_{j_0})=0$ and we are done. Therefore, $Z$ admits an inevitable biorthogonal system for every $0<\delta<1/2$. Now, to conclude that $Z$ is arbitrarily distortable we  just need to apply \cite[Proposition 13.28]{BenLin}.
\end{proof}
\section{Examples}
Let us show that a large class of twisted Hilbert spaces satisfy the Principle for Small Perturbations and thus we may apply Theorem \ref{main}. 
\begin{example}
The Kalton-Peck spaces $Z_p$ for $1<p<\infty$ admit an inevitable biorthogonal system for every $0<\delta<1/2$. In particular, they are arbitrarily distortable.
\end{example}
\begin{proof}
We only need to prove that $Z_p$ satisfies the P.S.P. and then apply Theorem \ref{main}. The idea that $Z_p$ satisfies the P.S.P. traces back to Kalton and Peck's paper \cite{KaPe}. The exact argument appears explicitly in the proof of \cite[Theorem 16.16.]{BenLin} of Benyamini and Lindenstrauss for $Z_2$. We reproduce it for the sake of completeness. First, we would like to recall an important fact that is: the sequence 
\begin{equation}\label{basis}
\{(e_1,0),(0,e_1),(e_2,0),...\}
\end{equation} is (in this order) a basis of $Z_p$, \cite[Theorem 4.10]{KaPe}. So, let $(w_j)_{j=1}^{\infty}$ be a sequence in $Z_p$ which is equivalent to the $\ell_p$ basis. Since $w_j\to 0$ weakly, we may apply The Bessaga-Pe\l czy\'nski Selection Principle \cite[Proposition 1.3.10]{KaPe} and pass to a subsequence (we do not relabel) so that there exists blocks $(x_j,u_j)$ of the basis (\ref{basis}) with $$\|w_j-(x_j,u_j)\|_{Z_p}\leq 2^{-j},\;\;j\in \mathbb N.$$ We claim that $\|u_j\|_{p}\to 0$. Let us observe that we always have the following inequality that will be crucial
\begin{eqnarray*}\label{crucial}
\left\| \left(\sum_{j=1}^n x_j,\sum_{j=1}^n u_j\right) \right\| &\geq& \left\| \sum_{j=1}^n x_j-\Omega_p\left( \sum_{j=1}^n u_j \right) \right\|\\
&\geq &  \left\| \Omega_p\left(\sum_{j=1}^n u_j\right)-\sum_{j=1}^n\Omega_p\left( u_j \right) \right\|-  \left\| \sum_{j=1}^n (x_j-\Omega_p\left(u_j\right) \right\|,
\end{eqnarray*}
where $\Omega_p$ denotes the Kalton-Peck map described in (\ref{kpmap}). The inequality above turns out to be impossible if $\inf_{j\in\mathbb N}\|u_j\| >0$. Indeed, $$\left\| \sum_{j=1}^n \left(x_j-\Omega_p(u_j)\right)\right\|^p_p\leq \sum_{j=1}^n \|x_j-\Omega_p(u_j) \|^p_p\leq n,$$ while 
\begin{equation}\label{equiva}
\left\| \left(\sum_{j=1}^n x_j,\sum_{j=1}^n u_j\right) \right\| \sim n^{\frac{1}{p}}
\end{equation} because the blocks $(x_j,u_j)$ are a small perturbation of $(w_j)$. If we assume that there is some $\eta>0$ so that $\|u_j\|\geq \eta >0$ for all $j\in \mathbb N$, then the blocks $u_j$ are semi-normalized and thus 
$$\left\| \Omega_p\left(\sum_{j=1}^n u_j\right)-\sum_{j=1}^n\Omega_p\left( u_j \right) \right\|\geq \rho n^{\frac{1}{p}}\log n,$$
by the definition of $\Omega_p$ for some $\rho>0$. If we insert these estimates in the crucial inequality, we find on the right side the lower bound $$\rho n^{\frac{1}{p}}\log n-n^{\frac{1}{p}},$$
which is impossible because the left side is of the order $n^{\frac{1}{p}}$. Thus, it must be $\|u_j\|_p\to 0$ as $j\to \infty$ and the claim is proved. So let us assume, passing to a subsequence if necessary, that $\|u_j\|_p\leq 2^{-j}$ and hence 
\begin{eqnarray*}
\|w_j-(x_j-\Omega_p(u_j),0)\|&\leq& \|w_j-(x_j,u_j)\|+\|(\Omega_p(u_j),u_j)\|\\
&\leq& 2^{-j}+\|u_j\|_p\\
&\leq& 2^{-j+1}.
\end{eqnarray*}
Thus, pick $y_j:=x_j-\Omega_p(u_j)$. Replacing $w_j$ by $\lambda_j w_j$, where the sequence $\lambda_j$ satisfies $\inf_j \lambda_j>0$ and $\sup_j \lambda_j<\infty$, we may assume that $\|y_j\|_p=1$. Therefore, $Z_p$ satisfies the P.S.P. and we are ready to apply Theorem \ref{main}.
\end{proof}
The rest of our examples deal with twisted sums of $\ell_2$ which are usually called twisted Hilbert spaces. Kalton provided us in \cite{ka} with a procedure so that for a given space $X$ with a shrinking unconditional basis, we may construct a twisted Hilbert space that we denote $Z(X)$. In its simplest form we could describe the process as follows: First, the complex interpolation formula $(X,X^*)_{1/2}=\ell_2$ holds true, see \cite{Wat} or \cite{CoS}. This means that every element $x\in \ell_2$, say with $\|x\|_2=1$, may be represented by certain holomorphic function $F$; here ``represented" means that $F(1/2)=x$ and $\|F\|\leq 2$ (this last norm is taken on the so-called Calder\'on space where the norm of $X$ is relevant). It is easy to check that the assignment $$\Omega(x):=F'(1/2)$$ is a centralizer. Thus, we may select our pairs as before and end up with $Z(X)$. It is not surprising in view of the formula above that $Z(X)$ is usually termed as the derived space of $X$. It is also not hard to check that $Z(\ell_r)=Z_2$ for $r\neq 2$, and $Z(c_0)=Z_2$.  The survey of Kalton and Montgomery-Smith \cite{KM} gives a short but accurate version of the process. 

We would like to say that unlike the Kalton-Peck map, we have no an explicit description of such $\Omega$ in general, we only know the existence of such map by the result of Nigel quoted above. This fact complicates very much the computation of the norm of an element $\|(x,y)\|$ in a twisted Hilbert space $Z$ with centralizer $\Omega$. However, using recent results of \cite{JS5}, we have a large class of twisted Hilbert spaces for which it is `easy' now to construct an inevitable biorthogonal system. This class contains the derived spaces of Maurey's class of asymptotic $\ell_p$-spaces \cite{Mau2}. Recall that a space $X$ with an unconditional basis is an asymptotic $\ell_p$-space if there is $C>0$ such that for every $n\in\mathbb N$ and every $n$-normalized blocks of the basis with $n<u_1<...<u_n$, we have that $[u_j]_{j=1}^n$ is $C$-equivalent to the natural basis of $\ell_p^n$.
\begin{example}Let $X$ be a space with a shrinking unconditional basis $(e_j)_{j=1}^{\infty}$. Assume $X$ satisfies one of the two following conditions:
\begin{enumerate}
\item $X$ is an asymptotic $\ell_r$-space with $1\leq r\neq 2<\infty$.
\item The only spreading model of a block basic sequence of $(e_j)_{j=1}^{\infty}$ is $c_0$.
\end{enumerate} Then, the twisted Hilbert space $Z(X)$ admits an inevitable biorthogonal system for every $0<\delta<1/2$. In particular, it is arbitrarily distortable.
\end{example}
\begin{proof}
The idea of the proof is the following. Since $X$ behaves asymptotically as, say $\ell_r$, then $Z(X)$ behaves asymptotically as $Z(\ell_r)=Z_2$. This notion of limit allows to share information between the spaces which is given by arbitrary large but finite data segments. As it is witnessed by the proof of Example 1, this is enough to get the P.S.P. so if $Z_2$ has the P.S.P. then so has $Z(X)$. Formally speaking, for the first case, we apply \cite[Theorem 4]{JS5} and find that the only spreading model of $Z(X)$ is $Z_2$. For the second, we use \cite[Proposition 9]{JS5} and also find that the only spreading model of $Z(X)$ is $Z_2$. Anyway, we may apply now \cite[Theorem 5]{JS5} and deduce that either in case (1) or (2) the space $Z(X)$ satisfies the P.S.P., so we conclude using Theorem \ref{main}. 
\end{proof}
The result above does not deal with $r=2$ but we may still give a nontrivial example of this; the Hilbert copy $\ell_2$ is a trivial example of this phenomena since it is trivial to see that $Z(\ell_2)=\ell_2$. Let $T$ denotes the Tsirelson space and $T^2$ the $2$-convexification which is well known to be an asymptotic $\ell_2$-space with an unconditional basis (see \cite{CS}). In particular, by Maurey's result \cite{Mau2}, the space $T^2$ admits an inevitable biorthogonal system. The derived space $Z(T^2)$, see \cite{JS} for more details on this space, also does.
\begin{example}
The twisted Hilbert space $Z(T^2)$ admits an inevitable biorthogonal system for every $0<\delta<1/2$. In particular, it is arbitrarily distortable.
\end{example}
\begin{proof}
We only need to prove that $Z(T^2)$ satisfies the P.S.P. and then apply Theorem \ref{main}. The proof that $Z(T^2)$ satisfies the P.S.P., which is indeed quite technical, is to be found at \cite[Theorem 3]{JS5}. 
\end{proof}
%The space $Z(T^2)$ is an asymptotic $\ell_2$-space XXX with no unconditional basis which is arbitrarily distortable.

\end{document}